\documentclass[12pt]{amsart}
\usepackage{graphicx}
\usepackage{tikz}
\usepackage{amsmath}
\usepackage[margin=1in]{geometry}

\begin{document}
	
	\begin{abstract}
		In this article, we investigate how Euler might have been led to conjecture the Prime Number Theorem, based on what he knew. We also speculate on why he did not do so.
	\end{abstract}
	
	\title{Could Euler have conjectured the prime number theorem?}
	\author{Simon Rubinstein-Salzedo}
	\address{Euler Circle, Palo Alto, CA 94306}
	\email{simon@eulercircle.com}
	%\address{Euler Circle, Palo Alto, CA 94306}
	%\email{simon@eulercircle.com}
	%\author{Author's name to be added later}
	\date{\today}
	\maketitle
	
	Euler likely was not motivated to study the distribution of prime numbers. However, conjecturing the Prime Number Theorem wouldn't have been that hard, based on what he knew, in particular the approximate growth rate of the partial sums $\sum_{p\le x} \frac{1}{p}$, where $p$ runs over only the primes up to $x$. Let us take a look at what Euler did do, and how he could have pushed it a little further, had he been sufficiently interested.
	
	\section{The sum of reciprocals of primes}
	
	Euler showed in his paper ``Variae observationes circa series infinitas''~\cite[Theorema 19]{euler1737} (Euler Archive number E72), presented to the St.\ Petersburg Academy in 1737 and published in 1744, that the sum of the reciprocals of the primes diverges, and he even worked out the growth rate of $\sum_{p\le x} \frac{1}{p}$. As was typical of Euler, his method of discovering this was beautiful, even if his lack of respect for divergent series may look horrifying to a modern mathematician. We present something similar to what Euler did in the remainder of this section; we have taken logarithms of every expression in his paper, to be more consistent with modern presentations.
	
\begin{figure}
		%\begin{center}\includegraphics{mascheroni.pdf}\end{center}
		\begin{tikzpicture}
		\draw (-.5,0) -- (8.5,0);
		\draw (0,-.5) -- (0,6);
		\filldraw (1,0) rectangle +(1,5);
		\filldraw (2,0) rectangle +(1,5/2);
		\filldraw (3,0) rectangle +(1,5/3);
		\filldraw (4,0) rectangle +(1,5/4);
		\filldraw (5,0) rectangle +(1,5/5);
		\filldraw (6,0) rectangle +(1,5/6);
		\filldraw (7,0) rectangle +(1,5/7);
		\filldraw[fill=white] (1,0) -- plot[variable=\t,samples=1000,domain=1:8] ({\t},{5/(\t)}) -- (8,0) -- (1,-0);
		\draw plot[variable=\t,samples=1000,domain=.9:8.5] ({\t},{5/(\t)});
		\draw (2,5) -- (2,0);
		\draw (3,2.5) -- (3,0);
		\draw (4,5/3) -- (4,0);
		\draw (5,5/4) -- (5,0);
		\draw (6,1) -- (6,0);
		\draw (7,5/6) -- (7,0);
		\foreach \n in {1,2,3,4}{
			\path node at (0,5*\n/4) {$-$};}
		\path node at (-.2,5) {$1$};
		\path node at (-.2,3.75) {$\frac{3}{4}$};
		\path node at (-.2,2.5) {$\frac{1}{2}$};
		\path node at (-.2,1.25) {$\frac{1}{4}$};
		\foreach \n in {1,...,8} {
			\path node at (\n,-.2) {$\n$};}
		% \path node at (5,4) {$\displaystyle\gamma=\lim_{n\to\infty} \left(\sum_{k=1}^n \frac{1}{k}-\int_1^{n+1} \frac{dx}{x}\right)$};
		% \path node at (6,2.5) {\Large $0<\gamma<1$};
		\end{tikzpicture}
		\caption{A ``proof by picture'' showing that $0<\gamma<1$.}
		\label{fig:gamma}
	\end{figure}
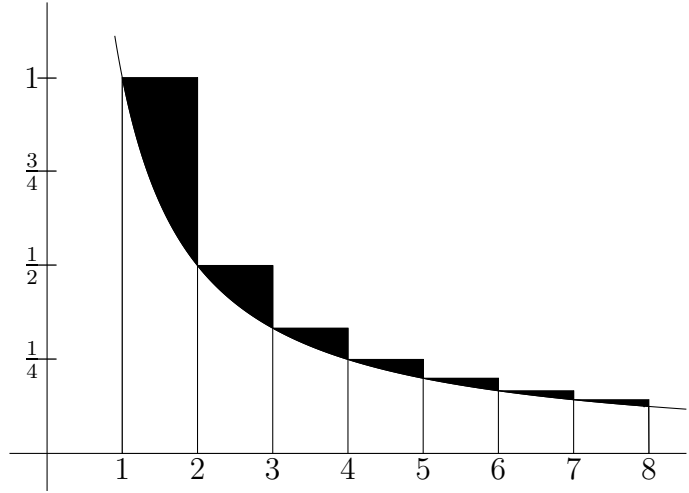

	It was already known that the harmonic series grows like $\log x$; that is, \[\sum_{n=1}^x \frac{1}{n}\approx \log(x),\] or, more precisely, \[\lim_{x\to\infty}\left(\sum_{n=1}^x\frac{1}{n}-\log(x)\right)=\gamma\approx 0.577.\] Figure~\ref{fig:gamma} shows that $0<\gamma<1$.
	
	Euler also knew a product formula, now known as the Euler product, for what would later be named the Riemann zeta function: \begin{equation}\sum_{n=1}^\infty \frac{1}{n^s}=\prod_{p\text{ prime}} \left(1-\frac{1}{p^s}\right)^{-1}\label{eulerproduct}\end{equation} whenever $\Re(s)>1$. This formula has a delightful interpretation as an analytic formulation of the fundamental theorem of arithmetic; this connection is described in many places, for instance~\cite{Apostol97}. Despite the fact that neither side of~(\ref{eulerproduct}) makes sense when $s=1$, Euler bravely set $s=1$ to obtain \begin{equation} 1+\frac{1}{2}+\frac{1}{3}+\frac{1}{4}+\cdots=\left(1-\frac{1}{2}\right)^{-1}\left(1-\frac{1}{3}\right)^{-1}\left(1-\frac{1}{5}\right)^{-1}\left(1-\frac{1}{7}\right)^{-1}\cdots. \label{eulerprods=1}\end{equation} Formula~(\ref{eulerprods=1}) is Theorema 7 in~\cite{euler1737}. Taking logarithms, he then obtained \begin{equation} \log\left(1+\frac{1}{2}+\frac{1}{3}+\frac{1}{4}+\cdots\right)=-\log\left(1-\frac{1}{2}\right)-\log\left(1-\frac{1}{3}\right)-\log\left(1-\frac{1}{5}\right)-\log\left(1-\frac{1}{7}\right)-\cdots.\label{logEP}\end{equation} Next, he used the power series expansion \[-\log(1-x)=x+\frac{x^2}{2}+\frac{x^3}{3}+\frac{x^4}{4}+\cdots\] and regrouped terms, so that the right side of~(\ref{logEP}) becomes \[\left(\frac{1}{2}+\frac{1}{3}+\frac{1}{5}+\frac{1}{7}+\cdots\right)+\frac{1}{2}\left(\frac{1}{2^2}+\frac{1}{3^2}+\frac{1}{5^2}+\frac{1}{7^2}+\cdots\right)+\frac{1}{3}\left(\frac{1}{2^3}+\frac{1}{3^3}+\frac{1}{5^3}+\frac{1}{7^3}+\cdots\right)+\cdots.\] Since all the groups but the first add up to something finite, they will not make a significant contribution in the limit, and Euler ignored them. Hence, in the limit, \[\frac{1}{2}+\frac{1}{3}+\frac{1}{5}+\frac{1}{7}+\cdots=\log\left(1+\frac{1}{2}+\frac{1}{3}+\frac{1}{4}+\cdots\right)=\log\log\infty.\] This result is Theorema 19 in~\cite{euler1737}. In modern notation, we would write \begin{equation} \sum_{\substack{p\le x\\ p\text{ prime}}} \frac{1}{p}\approx\log\log x. \label{growthrate}\end{equation} One can get rid of all the divergent series and odd-looking expressions like $\log\log\infty$ by taking limits and being careful with error terms, as a modern mathematician would do, but that wasn't a necessary part of the mathematical culture in the $18^\text{th}$ century.
	
	\section{On to the prime number theorem}
	
	At this point, Euler stopped. However, he could have used~(\ref{growthrate}) to estimate the number of primes up to $x$, at least conjecturally, using a quick and easy computation. Prime numbers were not the focus of his paper~\cite{euler1737}, but he could easily have returned to the topic at some later point, had he chosen to do so.
	
	Let us try to estimate the number of primes between $x$ and $kx$, for some $k>1$, and let us write $\pi(x)$ for the number of primes up to $x$. Based on~(\ref{growthrate}), we estimate \[\sum_{\substack{x<p\le kx\\ p\text{ prime}}}\frac{1}{p}\approx \log\log(kx)-\log\log(x).\] If $k$ is only slightly bigger than 1, then each term in the sum on the left is roughly $\frac{1}{x}$, and the number of terms in the sum is $\pi(kx)-\pi(x)$. Hence, we have \begin{equation}\log\log(kx)-\log\log(x)\approx\frac{\pi(kx)-\pi(x)}{x}.\label{loglog}\end{equation} The goal is now to estimate $\log\log(kx)-\log\log(x)$. Rewriting~(\ref{loglog}), we have \[\pi(kx)-\pi(x)\approx x(\log\log(kx)-\log\log(x))=x\log\left(1+\frac{\log k}{\log x}\right).\] Expanding the outer logarithm as a power series and truncating after the first term, we obtain \[\pi(kx)-\pi(x)\approx \frac{x\log k}{\log x}.\] Now, choosing $k=1+\frac{1}{x}$ and again taking the first term of the power series for $\log(k)=\log\left(1+\frac{1}{x}\right)\approx\frac{1}{x}$, we get \begin{equation} \pi(x+1)-\pi(x)\approx\frac{1}{\log(x)}.\label{probprime}\end{equation} Now, we can determine $\pi(x)$ by summing: \[\pi(x)=\sum_{n=1}^{x-1}(\pi(n+1)-\pi(n))\approx \sum_{n=2}^{x-1}\frac{1}{\log(n)}\approx\int_2^x \frac{dt}{\log t},\] which is one form of the prime number theorem. (We change the lower index of the sum and integral to 2 to avoid the problem of division by 0 that occurs when $n=1$. This does not affect the asymptotics.) A slightly different heuristic for the prime number theorem also coming from the sum of the reciprocals of the primes can be found in Sandifer's March 2006 column of ``How Euler Did It''~\cite{Sandifer06}. Yet another way of conjecturing the prime number theorem from the series comes from Abel's technique of partial summation, not yet developed in Euler's time, which allows one to compute $\sum_{n\le x} c_nf(n)$ from the sum $\sum_{n\le x} c_n$, here used with $c_n=\frac{1}{n}$ if $n$ is prime and 0 otherwise, and $f(n)=n$; see for instance~\cite[\S 1.3]{CM05}.
	
	\section{Primes and random variables}
	
	So, why didn't Euler do this? One can only speculate. Perhaps he was simply not sufficiently interested in the question of prime density; after all, as arguably the most prolific mathematician in history, he clearly had enough other things to occupy his attention! Indeed, his paper~\cite{euler1737} was not primarily about prime numbers, but rather about infinite series and their relation to infinite products, a subject that he was deeply interested in throughout his life, culminating in such gems as his solution to the Basel problem of evaluating the sum $\sum_{n=1}^\infty \frac{1}{n^2}$~\cite{euler1740}, E41, and the pentagonal number theorem~\cite{euler1783}, E541. Or perhaps he simply had no reason to believe that there would be any large-scale patterns in the primes. But my suspicion is that something like~(\ref{probprime}) would have felt very strange to an $18^\text{th}$-century mathematician. After all, what can it possibly mean? The left side is either 0 or 1: it's 0 if $x+1$ is composite and 1 if $x+1$ is prime. But the right side is some real number between 0 and 1, which does not ``know'' anything about prime numbers; rather, it's decaying smoothly. So~(\ref{probprime}) is nonsense, but, to quote Gilbert and Sullivan~\cite{gs2010}, ``oh, what precious nonsense!''
	
	Modern mathematicians understand things like~(\ref{probprime}) differently. There is no randomness in the primes: a fixed number $n$ is either prime or composite, and no amount of coin flipping or die rolling will ever change that. However, \emph{statistics} about primes behave like statistics about random numbers, with a certain distribution. Imagine building a set $S$ of positive integers as follows: for each integer $n\ge 2$, pick a random number $\alpha_n$ uniformly on $[0,1]$, and put $n$ into $S$ if $\alpha_n<\frac{1}{\log n}$, with each $\alpha_n$ being chosen independently from all the rest. This model is known as the Cram\'er model, after Harald Cram\'er, who wrote several papers such as~\cite{Cramer20} describing it. This model is elegantly described and analyzed in a recent survey article~\cite{Sound07} of Soundararajan.
	
	For instance, the Cram\'er model predicts the Twin Prime Conjecture, which is still open, in spite of recent progress of Goldston--Pintz--Y\i ld\i r\i m~\cite{GPY09}, Zhang~\cite{Zhang14}, Maynard~\cite{Maynard15}, and the Polymath project~\cite{Polymath14}. The Twin Prime Conjecture states that there are infinitely many positive integers $n$ so that $n-2$ and $n$ are both prime. The argument goes as follows: the probability that $n-2$ and $n$ are both in $S$ is $\frac{1}{\log(n-2)}\times\frac{1}{\log(n)}\ge\frac{1}{\log(n)^2}$. Thus, the expected number of pairs $n-2,n\in S$ is \[\sum_{n=4}^\infty \frac{1}{\log(n-2)\log(n)}\ge\sum_{n=4}^\infty \frac{1}{\log(n)^2}=\infty.\]
	
	This approach must be used with care: for example, the same argument as in the preceding paragraph can be used to show that $S$ is expected to contain infinitely many pairs of consecutive numbers, whereas 2 and 3 are the only consecutive primes, so statistics about pairs of consecutive elements of $S$ and pairs of consecutive primes behave completely differently. However, one might be tempted to say that any statistics about $S$ and the primes agree ``unless there is an obvious reason that they don't.'' Of course this is only a rough heuristic: after all, one person's obvious reason may be another person's deep result.
	
	While the Cram\'er model is only a heuristic, it is sometimes possible not only to make conjectures but also to prove theorems by taking a probabilistic approach to primes and divisibility. One of the first examples of these theorems was the Erd\H{o}s--Kac Theorem~\cite{EK40}, which shows that the number of prime factors of numbers of size roughly $n$ becomes normally distributed as $n\to\infty$, with mean and variance $\log\log(n)$. Since then, the probabilistic method has become a standard technique in analytic number theory, as well as many other areas of mathematics.
	
	\section*{Acknowledgments}
	
	I would like to thank the anonymous referees for helpful suggestions that improved the article, and in particular for bringing the article~\cite{Sandifer06} to my attention.
	
	\bibliographystyle{alpha}
	\bibliography{euler}

\end{document}